\overfullrule=0pt

\magnification 1200

\def \dd#1{{\bf#1}}

\def\cl#1{{\cal#1}}



\def\ouv#1{\smash{\mathop{#1}\limits^{\lower 1pt\hbox
{$\scriptscriptstyle\circ$}}}}

\def\hfl#1#2{\smash{\mathop{\hbox to 12mm{\rightarrowfill}}
\limits^{\scriptstyle#1}_{\scriptstyle#2}}}


\long\def\eno#1#2{\par\smallskip{\bf{#1}}{\it\ {#2}}\par\medskip}

\def\stit#1{\vskip 3mm plus 1mm minus 2mm {\bf{#1}}
		\smallskip}

\font\tir=cmbx10 at 12pt

\def\ref#1#2#3#4{{\bf #1}{\ #2}{\it ,\ #3}{,\ #4}\medskip}


\def \picture #1 by #2 (#3){\midinsert \centerline 
{\vbox to #2{\hrule width #1 heigth 0pt 
depth 0pt \null \vfill \special {picture #3}}}\endinsert }

\def\scaledpicture #1 by #2 (#3 scaled #4) {{
\dimen0 =#1 \dimen1 =$2
\divide \dimen0 by 1000 \multiply \dimen0 by #4
\divide \dimen1 by 1000 \multiply \dimen1 by #4
\picture \dimen0 by \dimen1 (#3 scaled $4)}}

\def\figure #1 #2 #3 {\midinsert \vglue 3mm 
{\vbox to #3 {\hrule width 6cm height 0cm depth 0cm \vfill
{\special {picture #1 scaled #2}}}}\vglue 2mm \endinsert}

\magnification=1200


\def\ouv#1{\smash{\mathop{#1}\limits^{\lower 1pt\hbox
{$\scriptscriptstyle\circ$}}}}


\long\def\eno#1#2{\par\smallskip{\bf{#1}}{\it\ {#2}}\par\medskip}

\def\stit#1{\vskip 3mm plus 1mm minus 2mm {\bf{#1}}
		\smallskip}

\def\TTtit#1{\centerline{\tir {#1}}}

\font\tir=cmbx10 at 12pt

\def\ref#1#2#3#4{{\bf #1}{\ #2}{\it ,\ #3}{,\ #4}\medskip}


\vskip 2cm

\TTtit {Log-euclidean geometry and }
\TTtit {"Grundlagen der Geometrie"} 
\vskip .7cm plus 2mm minus 2mm

{\centerline { {\bf Ricardo P\'erez-Marco*
}}}

\vskip .5cm plus 1mm minus 1mm

{

{\bf Abstract.} We define the simplest log-euclidean geometry. This geometry
exposes a difficulty hidden in Hilbert's list of axioms presented
in his "Grundlagen der Geometrie". The list of axioms appears to be
incomplete if the foundations of geometry are to be independent of
set theory, as Hilbert intended. In that case we need to add a missing axiom.
Log-euclidean geometry satisfies
all axioms but the missing one, the fifth axiom of congruence
and Euclid's axiom of parallels. This gives an elementary proof (with no
need of Riemannian geometry) of the independence of these axioms from
the others.
\footnote{}{*CNRS, IMJ-PRG, E-mail: ricardo.perez-marco@imj-prg.fr. This manuscript was composed in 2004 and two references are updated.}
\footnote{}{2010 Mathematics Subject Classification 2010: 51M05, 30F99, 03B30, 01A20.}
\footnote{}{Key Words: Euclidean geometry, Euclid's axioms, Grundagen, log-Riemann Surfaces.}

%
%
%
}

\medskip

\stit {1. Log-euclidean geometry.}

\stit {1.1. A metric space.}

Consider two pointed euclidean planes, that we can identify with the complex
plane $\dd C-\{ 0 \}$, slitted along an euclidean half-line, say along
the negative real axes $]-\infty, 0[$. We glue by the identity on charts
the upper, resp. lower,
boundary of the slit of one plane with the lower, resp. lower, boundary of the
slit in the other plane. We add the missing point at $0$, that
we continue to denote by $0$.
The distance between two points in this space is the infimum of lengths of
paths joining these two points. The length is being computed in each euclidean
chart. We denote by $\cl X$ the metric space constructed. Conformally
the concrete Riemann surface\footnote{$\dagger$}{A concrete Riemann surface is
a Riemann surface endowed with a canonical system of charts.}
obtained is the Riemann surface associated to
$z\mapsto \sqrt {z}$. This is why we refer to $0$ as the ramification point.
Note that it is necessary to add a the ramification point $0$ in order to 
have a complete space. It is elementary to prove:

\eno {Proposition.}{Given two distinct points $A,B \in \cl X$ there is
a unique geodesic path joining them. This geodesic path is
an euclidean segment, or a broken line formed by two euclidean segments
meeting at the ramification point $0$.}

\stit {1.2. Log-euclidean points and lines.}

In log-euclidean geometry a point is an element of the space $\cl X$.
Lines in log-euclidean geometry are defined to be euclidean lines or two half-lines
meeting at the ramification point $0$ (this point being included).

\medskip

There exists a log-euclidean geometry in each log-Riemann surface.
These are Riemann surfaces with a set of canonical charts
of a certain type.
Log-euclidean metric characterizes the log-Riemann surface structure
(see [Bi-PM1] and [Bi-PM2] for more information on log-Riemann surfaces).

\stit {2. Hilbert axioms.}

\stit {2.1) Preliminary comments.}

The purpose of the "Grundlagen der Geometry" is to give a complete
list of independent axioms upon which euclidean geometry can be
build. This philosophy is more than 2300 years old and is
borrowed from "The Elements" of Euclid ([Eu]) which marks the birth of
formal Mathematics. The purpose of
the "Grundlagen" is to add those axioms
that were assumed to be self-evident in "The Elements" and to prove
the independence of the axioms. The problem of the independence
of Euclid's fifth postulate on parallels has kept busy generations
of mathematicians (Saccheri, Lambert, Legendre,...) until the creation
of non-euclidean geometries by Gauss, Bolyai, Lobatchevski and their
analytic realization. We refer to [Bo] and [Ef] for a delightful historical
survey.

\medskip

Hilbert's Grundlagen sets the axioms for three dimensional euclidean
geometry. A subset of these axioms (named planar axioms) compose 
the axiomatic system of two dimensional euclidean geometry. We restrict to
discuss planar
geometry in this article.

\medskip

As Hilbert states in his first sentence, he considers in planar
geometry two types of objects:
{\it points}, denoted with capital letters $A, B, C, \ldots$, and {\it lines},
denoted by lower case letters $a, b, c,\ldots$ The goal of the axioms is to
describe relations between these objects. Thus, with the system of axioms
at hand, a blind person can prove theorems in euclidean geometry.  An
important subtle point of Hilbert's conception is that the nature of
"points" or "lines" is irrelevant. A famous quote of Hilbert
that report several sources ([We], [Re], [BB], see also [Hi-Ce]) is the
following

\bigskip

{\it "...One should be able to talk about chairs, tables and mugs of beer
instead of points, lines or planes..."}

\bigskip

All of this discusion shows that in Hilbert's system
{\bf lines are not considered to be sets of points} to start with
(as Hilbert would say "tables are not sets of chairs".)
There is a very good reason for this point of view.
Hilbert carefully avoids to base axiomatic euclidean geometry
on axiomatic set theory. He is well aware (since the work of Cantor)
that axiomatic set theory is far more
complex to build.

\stit {2.2) Hilbert's axioms.}

Hilbert's axioms are divided into five distinct groups:

\medskip

{\parindent=2cm
\item {I.} Axioms of connection.
\item {II.} Axioms of order.
\item {III.} Axioms of congruence.
\item {IV.} Axiom of parallels.
\item {V.} Axioms of continuity.
}
\medskip

There have been ten editions of the Grundlagen (see [Hi-Ro]), seven
during Hilbert's lifetime. The set of axioms evolved to reach the
final form after the seventh edition. In particular in the second
edition one axiom became a theorem proved by E. H. Moore [Mo].

We first focuss on the first group of axioms of connection in
planar geometry.

\stit {2.3) Log-euclidean geometry and the axioms of connexion.}

\stit {2.3.1) Axioms of connection.}

In the words of Hilbert, the axioms of connection build
connections or links between point objects and line objects.

The first two planar axioms are the following:

\medskip

{\bf I-1. Given two points $A$ and $B$, there always exists a line
$a$ that corresponds to the two points $A$ and $B$.}

\medskip

{\bf I-2. Given two points $A$ and $B$, there is only one line that
corresponds to the two points $A$ and $B$.}

\medskip

Hilbert takes good care in indicating that "two points" means
"two distinct points"; and that "corresponds" will be or can be
replaced by
"goes through", or "$A$ is located in $a$", or "$A$ is a point of $a$",
or "$A$ belongs to $a$", etc

\medskip

As explained above "$A$ belongs to $a$" cannot be
taken in a set theoretical sense. This seems to have caused
some psychological difficulty: According to [Hi-Ro] after the sixth
edition of the Grundlagen
"belongs to" was replacing "corresponds to" from the previous edition.
And this term is translated in French as "appartient". But in [Hi-Ce],
which is a translation of the seventh edition, it is translated in
Spanish as "se corresponde". The German word used in the text
is "zusammengeh\"oren".

But what word to use, or in its place to use a potato to mark the connection,
is after all just a matter of semantics and is irrelevant.

\stit {2.3.2) Log-euclidean geometry.}

We have to define the correspondence stated in axioms {\bf I-1}
and {\bf I-2}
in log-euclidean geometry. The line corresponding to two points
$A$ and $B$ will be the euclidean line passing through these points
if an euclidean segment is the geodesic joining these two points (including
the case when this euclidean segment contains the ramification point
$0$). In the case when a broken segment is the geodesic joining these
two points then the corresponding line is obtained by extending the two
segments in the oposite direction of the ramification point $0$.

\medskip

It is clear that such correspondence satisfies the axioms {\bf I-1}
and {\bf I-2}.
Moreover the third planar connection axiom is also trivially verified:

\medskip

{\bf I-3. Any line goes through at least two points. There are at least
three points not in a line.}

\medskip

Again here the "goes through" and "belonging" statement are to be taken
in the sense specified by Hilbert's remarks. Thus, for example,
to say that a line
goes through at least two points $A$ and $A'$ means that there are points
$B$, distinct from $A$, and $B'$, distinct from $A'$, such that the
line is associated to $A$ and $B$, and also to
$A'$ and $B'$.

\stit {2.3.3) Log-euclidean geometry and the planar part of the
first Theorem  of the Grundlagen.}

After stating the axioms of this first group, Hilbert states without
proof his first Theorem. The first part of the theorem is a planar
statement that supposedly follows from the planar axioms:

\eno{Theorem 1 (Planar part).}{Two lines have one or no points in common.}

In an equivalent form, if two lines have two distinct
points in common then they must be equal.

Also in [Ef] p.43 this part of Theorem 1 is stated to follow trivially 
from axiom {\bf I-2} and does not deserve a proof.

\medskip

Now consider two points $A$ and $B$ in log-euclidean geometry joined
by a geodesic which is a broken segment. The line $a$ associated
to $A$ and $B$, and the line $b$ associated to $A$ and the ramification
point $0$, do have a half line in common!

\medskip

{\bf Conclusion: Theorem 1 does not hold in log-euclidean geometry, thus
theorem 1 cannot follow from axioms I-1, I-2 and I-3.}

\stit {3) Solution of the paradox.}

The confusion in Hilbert's text is created by the terminology "belonging
to a line" that is not to be taken in a set theoretic sense, but apparently
is taken in this sense in order to prove Theorem 1.
If we assume set theory, and we formulate the second
axiom as

\medskip

{\bf I-2'. Given two points $A$ and $B$ there exists at most one line
containing $A$ and $B$.}

\medskip

where we understand that "containing" is taken in a set theoretic sense
(thus this implies that the objects "lines" are sets containing
some "point" objects), then the proof of theorem 1 is straithforward.

\medskip

This appears to be a subtle point, but log-euclidean geometry
shows that one must be careful on this point. 

\medskip

One needs to add a missing axiom in order to exclude log-euclidean
geometry. We propose:

\medskip

{\bf I-2bis. If a line $a$ corresponds to $A$ and $C$, and also
to $B$ and $D$, then if $A$ and $B$ are distinct points, the line $a$
also corresponds to $A$ and $B$. }

\medskip

Note that axiom {\bf I-2bis} is not fulfilled by log-euclidean geometry.
Just consider a broken line with $A$ and $B$ in one of the half lines
determined by the ramification point $0$ and and $C$ and $D$ in the other
half line.

We can now prove theorem 1.

\stit {Proof of Theorem 1.}

Assume that the lines $a$ and $a'$ have two points $A$ and $B$ in common.
This means that there exist points $C$ and $D$ such that the line
$a$ corresponds to $A$ and $C$, and it also corresponds to $B$ and $D$.
Thus using the axiom {\bf I-2bis} it also corresponds to $A$ and $B$.  With
this same argument the line $a'$ also corresponds to $A$ and $B$.
Thus we conclude
that the lines $a$ and $a'$ coincide using axiom {\bf I-2}.$\diamond$

\medskip

It is curious to note that the weaker statement:

\medskip

{\bf "If the line $a$
corresponds to $A$ and $B$, and also to $A$ and $C$, and if $B$ is distinct
from $C$ then $a$ corresponds to $B$ and $C$"}
\medskip

apparently was included in axiom {\bf I-2} of the first edition of the Grundlagen
according to [Hi-Ro]. This still seems insufficient to prove Theorem 1
(but is enough to exclude log-euclidean geometry).

\medskip

In conclusion we want to insist that the subtle point raised here is
not irrelevant. If one aims to build an axiomatic system of euclidean
geometry not relying on axiomatic set theory (as Hilbert intended),
one should add
an axiom as {\bf I-2bis}.

\stit {4) Log-euclidean geometry and the other axioms.}

It is straightforward to check that the axioms of order (group II)
and continuity (group V)
are satisfied by log-euclidean geometry.
Also, taking as lenght of segments in log-euclidean geometry
the euclidean lengths, then the axioms of congruence relative
to segments (these are {\bf III-1}, {\bf III-2}, {\bf III-3})
are fulfilled.

\medskip

In order to check the axioms of congruence of angles {\bf III-4}
and {\bf III-5}
we need to define the magnitude of angles in log-euclidean geometry.
Lines meeting at a point distinct from $0$ form an angle whose magnitude
is the one given by euclidean geometry. For angles having the ramification
point $0$ as vertex, we define its magnitude as the euclidean angle formed
by the angle formed by the bissectors of the broken lines at $0$ pointing
towards the region where a full sheet is attached. With this convention
axiom {\bf III-4} is fulfilled.

\medskip

But axiom {\bf III-5} on congruence of triangles is not fulfilled as well
as Euclid's axiom of parallels {\bf IV-1}. Any line has infinitely many
parallel lines going through an external point.

\bigskip

{\bf Bibliography.}

\medskip

\ref{[Bi-Be]}{G. BIRKHOFF, M.K. BENNET}{Hilbert's Grundlagen der Geometrie}
{Rendiconti del Circolo Matematico di Palermo, {\bf 36}, II, p.343-389, 1987.}

\ref{[Bi-PM1]}{K. BISWAS, R. P\'EREZ-MARCO}{Log-Riemann surfaces}
{ArXiv 1512.03776, 2015.}

\ref{[Bi-PM2]}{K. BISWAS, R. P\'EREZ-MARCO} {Uniformization of simply connected finite type log-Riemann surfaces} {Geometry, groups and dynamics, 205-216, 
Contemp. Math., {\bf 639}, Amer. Math. Soc., Providence, RI, p.205-216, 2015.}

\ref{[Bo]}{R. BONOLA}{non-euclidean geometry}
{With a supplement containing "The theory of parallels" by N. Lobatchevski,
and "The science of absolute space" by J. Bolyai, Dover, New York, 1955.}

\ref{[Ef]}{N. EFIMOV}{G\'eom\'etrie sup\'erieure}
{French translation by E. MAKHO, Mir Editions, Moscow, 1985.}

\ref{[Eu]}{EUCLID}{The Elements}
{Critical edition J.L. Heiberg, Leipzig, Teubner, 1883; and translation 
by Heath, Camb. Univ. Press, 1908.}

\ref{[Hi]}{D. HILBERT}{Grundlagen der Geometrie}
{10 Editions: 1899, 1903, 1909, 1913, 1922, 1923, 1930, 1956, 1962, 1968.}

\ref{[Hi-Ce]}{D. HILBERT}{Fundamentos de geometr\'\i a}
{Spanish translation of the 7th german edition by Francisco 
Cebri\'an, with a historical introduction by J.M. S\'anchez Ron, 
Colecci\'on Textos Universitarios, {\bf 5}, CSIC, Madrid, 1991.}

\ref{[Hi-Ro]}{D. HILBERT}{Les fondements de la G\'eom\'etrie}
{Commented translation by Paul Rossier, Dunod, Paris, 1971.}

\ref{[Mo]}{E.H. MOORE}{}{Transactions of the Amer. Math. Soc., 1902.}

\ref{[Re]}{C. REID}{Hilbert}{Berlin-Heidelberg-New York, p.57, 1970.}

\ref{[We]}{H. WEYL}{David Hilbert and his mathematical work}
{Bull. Amer. Math. Soc., {\bf 50}, p.612-654, 1944.}

\bigskip

%
%
%

\end